\documentclass[10pt]{article}
\pdfoutput=1
\usepackage[backend=biber,
            bibstyle=phys,
            sorting=none,
            url=false,
            doi=false,
            eprint=true,
            sortcites=true,
            citestyle=numeric-comp]{biblatex}

\DeclareFieldFormat*{title}{\mkbibitalic{#1}}
\DeclareFieldFormat{labelnumberwidth}{\mkbibbrackets{#1}~~}

\usepackage[utf8]{inputenc}
\usepackage[top=1in,bottom=1in,left=1in,right=1in]{geometry}
\usepackage{authblk}
\usepackage{amsmath,amsfonts,amssymb,amsthm,xcolor}
\usepackage[breaklinks]{hyperref}

\newtheorem{proposition}{Proposition}

\numberwithin{equation}{section}

\newcommand{\R}{\mathbb{R}}

\newcommand{\Z}{\mathcal{Z}}
\newcommand{\X}{\mathcal{X}}
\newcommand{\A}{\mathcal{A}}
\newcommand{\Q}{\mathcal{Q}}

\newcommand{\D}{\Delta}

\newcommand{\ket}[1]{\left\rvert #1 \right\rangle}
\newcommand{\bra}[1]{\left\langle #1 \right\rvert}

\renewcommand*{\Affilfont}{\normalsize\small}
\author[1]{André Beaudoin\,}
\author[2]{Geoffroy Bergeron\,}
\author[3]{Antoine Brillant\,}
\author[4]{Julien Gaboriaud\,}
\author[5]{Luc Vinet\,}
\author[6]{Alexei Zhedanov\!
\vspace{.5em}}
\affil[1]{Département de génie physique, Polytechnique Montréal, Montréal (Québec), H3T
1J4, Canada.
\vspace{.9em}}
\affil[2,4,5]{Centre de Recherches Mathématiques, Université de Montréal,\protect\\
P.O. Box 6128, Centre-ville Station, Montréal (Québec), H3C 3J7, Canada.
\vspace{.9em}}
\affil[2,3,4,5]{Département de Physique, Université de Montréal, Montréal
(Québec), H3C 3J7, Canada.
\vspace{.9em}}
\affil[6]{School of Mathematics, Renmin University of China, Beijing, 100872, China.
\vspace{1.5em}
}
{
 \makeatletter
 \renewcommand\AB@affilsepx{: \protect\Affilfont}
 \makeatother
 \affil[ ]{E-mail addresses}
 \makeatletter
 \renewcommand\AB@affilsepx{, \protect\Affilfont}
 \makeatother
 \affil[1]{andre.beaudoin@polymtl.ca}
 \affil[2]{geoffroy.bergeron@umontreal.ca}
 \affil[3]{antoine.brillant@umontreal.ca}
 \affil[4]{julien.gaboriaud@umontreal.ca}
 \affil[5]{vinet@crm.umontreal.ca}
 \affil[6]{zhedanov@yahoo.com}
}

\title{Orthogonal polynomials and the deformed Jordan plane}
\date{\today}

\begin{document}
\maketitle
\thispagestyle{empty}
\hrule
\begin{abstract}
We consider the unital associative algebra $\mathcal{A}$ with two generators
$\mathcal{X}$, $\mathcal{Z}$ obeying the defining relation
$[\mathcal{Z},\mathcal{X}]=\mathcal{Z}^2+\Delta$.  We construct irreducible tridiagonal
representations of $\mathcal{A}$.  Depending on the value of the parameter $\Delta$, these
representations are associated to the Jacobi matrices of the para-Krawtchouk, continuous
Hahn, Hahn or Jacobi polynomials.
\end{abstract}
\centerline{\textbf{Keywords:} Para-Krawtchouk polynomials, deformed Jordan plane,
tridiagonal representations.}
\centerline{\textbf{MSC Class:} 33C47, 33C45, 39A13, 47L55.}
\centerline{}
\hrule

\section{Introduction}

This paper is devoted to the study of irreducible tridiagonal representations of the
two-generated algebra $\A$ which is a deformation of the Jordan plane.
% algebra $[\Z,\X]=\Z^2+\D$ with $\D$ a parameter.
It is shown how the para-Krawtchouk
polynomials appear quite naturally in this context, along with the other families of
classical orthogonal polynomials (OPs) of the Jacobi, continuous Hahn and Hahn type.

The algebra $\A$ over $\R$, with generators $\X$, $\Z$ and satisfying
\begin{align}\label{eq:adef}
 [\Z,\X]=\Z^2+\D
\end{align}
with $\D$ a parameter, is a special case of the most general two-generated quadratic
algebra $\Q$ with defining relation
\begin{align}\label{eq:genquad}
\alpha_1\X^{2}+\alpha_{2}\X\Z+\alpha_{3}\Z\X+\alpha_{4}\Z^{2}
+\alpha_{5}\X+\alpha_{6}\Z+\alpha_{7}=0.
\end{align}
This algebra has been of interest to various communities.  Ring theorists have provided
classifications \cite{Smith1992, Gaddis2015} of the special cases it entails and studied
their properties.  The algebra $\Q$ has also been related to non-commutative
probability theory \cite{BozejkoKummereretal1997} and is related to martingale polynomials
associated to quadratic harnesses \cite{BrycMatysiaketal2007}. On the physics side, $\Q$
describes various $1$D asymmetric exclusion models \cite{DerridaEvansetal1993,
EsslerRittenberg1996, UchiyamaSasamotoetal2004}.

Recently, the last two authors have begun connecting $\Q$ and its various isomorphism
classes to families of special functions. In \cite{TsujimotoVinetetal2017}, by studying
tridiagonal representations of the $q$-oscillator algebra $\X\Z-q\Z\X=1$,
% \RED{
they have identified how they encompass the recurrence relations of the big $q$-Jacobi,
the $q$-Hahn and the $q$-para-Krawtchouk polynomials.
% }
The case of the $q$-Weyl algebra
$\X\Z-q\Z\X=0$ has also been studied in \cite{Zhedanov1993}.  The present paper will add
to this program by considering an interesting special case of \eqref{eq:genquad} and
identifying the orthogonal polynomials that can be interpreted from this algebra.

Since their introduction in \cite{VinetZhedanov2012}, para-polynomials have been the
object of growing interest. Four families have been defined and studied offering
para-versions of the polynomials of Krawtchouk, $q$-Krawtchouk, Racah and $q$-Racah type.
While they do not fall in the category of classical orthogonal polynomials \footnote{They
obey a three term recurrence relation but a higher order difference equation.}, they are
understood as non-standard truncations of infinite-dimensional families of classical OPs
\cite{GenestVinetetal2013c, LemayVinetetal2016, LemayVinetetal2018}.  In addition to their
natural occurence in the study of perfect state transfer and fractional revival in quantum
spin chains \cite{VinetZhedanov2012, GenestVinetetal2016a, LemayVinetetal2016a}, recent
advances have identified these para-polynomials as the basis for finite-dimensional
representations of degenerations of the Sklyanin algebra \cite{GaboriaudTsujimotoetal2020,
BergeronGaboriaudetal2021, BergeronGaboriaudetal2021a}. They have also appeared in the
study of the Dunkl oscillator in the plane \cite{GenestIsmailetal2014}.
% \RED{
The main goal of this paper is to show that these para-Krawtchouk polynomials as well as
the Jacobi, continuous Hahn and Hahn polynomials arise in representations of the
two-generated algebra $\A$.
 % corresponding to a deformation of the Jordan plane.
% }

% Let $\A$ be the algebra over $\R$ generated by $\X$, $\Z$ with defining relation
% \begin{align}\label{eq:adef}
%  \Z\X-\X\Z-\Z^2-\D=0, \qquad \D\in\R.
% \end{align}
When $\D\neq0$, $\A$ as defined in \eqref{eq:adef} is a deformation of the Jordan plane
(with $\X$ and $\Z$ viewed as noncommutative coordinates).  Three cases will be
distinguished depending on whether $\D=0$, $\D>0$ or $\D<0$.  These three cases will be
studied separately and provide a complete picture of the connection between the algebra
\eqref{eq:adef} and orthogonal polynomials.

The presentation is organized as follows. Section \ref{sec:tridiag} will introduce the
tridiagonal representations of the algebra $\A$ and the non-degeneracy condition.
Standardized versions of $\A$ corresponding to $\D=0$, $\D<0$, $\D>0$ will then be
examined in the following sections.  The case $\D=0$ will be studied in section
\ref{sec:c0} and the Jacobi OPs will appear, while the case $\D>0$ and the continuous Hahn
polynomials will be the object of section \ref{sec:cpu}. Section \ref{sec:cmu} will focus
on the case $\D<0$ and will feature both the Hahn and the para-Krawtchouk polynomials.
Some concluding remarks and perspectives will close the paper.

\section{Tridiagonal representations of the algebra $\A$}\label{sec:tridiag}
Consider a tridiagonal representation of $\A$ where $\X\mapsto X$ and
$\Z\mapsto Z$. The actions of $X$, $Z$ on a semi-infinite orthonormal basis
$\ket{n}, n=0,1,2,\dots$ are taken to be of the form
\begin{subequations}\label{eq:actions}
\begin{align}
 X\ket{n} &= c_n \ket{n-1} + b_n \ket{n} + a_n \ket{n+1},\\
 Z\ket{n} &= u_n \ket{n-1} + v_n \ket{n} + w_n \ket{n+1},
\end{align}
\end{subequations}
with $c_0=u_0=0$.
To ensure that such a representation is irreducible we shall assume that the off-diagonal
coefficients are non-zero for $n>0$. Acting with \eqref{eq:adef} on the basis $\ket{n}$
and using the above definitions, one obtains
\begin{multline}
 (ZX-XZ-Z^2-\D)\ket{n}=(c_n u_{n-1}-c_{n-1} u_n-u_{n-1} u_n)\ket{n-2}\\
 +(b_n u_n - b_{n-1} u_n + c_n v_{n-1} - u_n v_{n-1} - c_n v_n - u_n v_n)\ket{n-1}\\
 +(-\D-a_{n-1} u_n + a_n u_{n+1} - v_n^2 + c_n w_{n-1} - u_n w_{n-1} - c_{n+1} w_n
 - u_{n+1} w_n)\ket{n}\\
 +(a_n v_{n+1} - a_n v_n + b_n w_n - b_{n+1} w_n - v_n w_n - v_{n+1} w_n)\ket{n+1}\\
 +(a_n w_{n+1} - a_{n+1} w_n - w_n w_{n+1})\ket{n+2}.
\end{multline}
For the actions in \eqref{eq:actions} to define a representation of $\A$, each side of the
above equation must vanish. As the basis vectors are orthonormal, one obtains the
following conditions on the coefficients of \eqref{eq:actions} that define the
representations:
\begin{align}
 0 &= c_n u_{n-1}-c_{n-1} u_n - u_{n-1} u_n,\label{eq:5}\\
 0 &= b_n u_n - b_{n-1} u_n + c_n v_{n-1} - u_n v_{n-1} - c_n v_n - u_n v_n,\label{eq:6}\\
 0 &= -\D-a_{n-1} u_n + a_n u_{n+1} - v_n^2 + c_n w_{n-1} - u_n w_{n-1} - c_{n+1} w_n
      - u_{n+1} w_n,\label{eq:7}\\
 0 &= a_n v_{n+1} - a_n v_n + b_n w_n - b_{n+1} w_n - v_n w_n - v_{n+1} w_n,\label{eq:8}\\
 0 &= a_n w_{n+1} - a_{n+1} w_n - w_n w_{n+1}.\label{eq:9}
\end{align}

\subsection{General solutions to the recurrence relations}
We now determine the general solutions to the above system of recurrence equations.
Dividing
\eqref{eq:5} by $u_n u_{n-1}$, one obtains
\begin{align*}
 \frac{c_n}{u_n}-\frac{c_{n-1}}{u_{n-1}}=1.
\end{align*}
This implies
\begin{align}\label{eq:solphi}
 \phi_n = \phi_0 + n, \qquad \phi_n\equiv\frac{c_n}{u_n}.
\end{align}
Equation \eqref{eq:9} can be solved similarly. Dividing by $w_n w_{n+1}$, one has
\begin{align}\label{eq:soldelta}
 \delta_n = \delta_0 - n, \qquad \delta_n \equiv \frac{a_n}{w_n}.
\end{align}
Rewriting \eqref{eq:6} and \eqref{eq:8} in terms of $\phi_n$ and $\delta_n$ and dividing
by $u_n$ or $w_n$, respectively, one obtains
\begin{align}
 b_{n-1}-b_n &= (\phi_n - 1) v_{n-1} - (\phi_{n} + 1) v_n,\label{eq:b1}\\
 b_{n+1}-b_n &= (\delta_n - 1) v_{n+1} - (\delta_n+1) v_n \label{eq:b2}.
\end{align}
To solve for $v_n$, shift the index of \eqref{eq:b1} and add \eqref{eq:b2} to find
\begin{align}\label{eq:vnnp1}
 0 = (\delta_n - \phi_{n+1} - 2)v_{n+1} - (\delta_n - \phi_{n+1} + 2)v_n.
\end{align}
Substituting the solutions \eqref{eq:solphi} and \eqref{eq:soldelta} in \eqref{eq:vnnp1}
leads to
\begin{align}
 0 &= (\delta_0 - \phi_0 - 2(n+2) +1)v_{n+1} - (\delta_0 -\phi_0 - 2n +1)v_n\\
  &= \mu_{n+2}v_{n+1} - \mu_n v_n,
\end{align}
with $\mu_n \equiv (\delta_0 -\phi_0 - 2n +1)$.
Multiplying the above by $\mu_{n+1}$ as an integrating factor, one can solve the
recurrence to obtain
\begin{align}\label{eq:solv}
 v_n = \frac{(\delta_0 -\phi_0 - 1)(\delta_0 -\phi_0 + 1) v_0}
            {(\delta_0 -\phi_0 - 2n +1)(\delta_0 -\phi_0 - 2n -1)}.
\end{align}
To find $b_n$, substract instead \eqref{eq:b1} with shifted index from
\eqref{eq:b2} and get
\begin{align}
 b_{n+1} - b_n = \frac{1}{2}(\delta_n + \phi_n + 1) (v_{n+1} - v_n),
\end{align}
which, upon using \eqref{eq:solphi} and \eqref{eq:soldelta}, can be solved immediately and
yields
\begin{align}\label{eq:solb}
 b_n %&=\frac{1}{2}(\delta_{n-1}+\phi_{n-1} + 1)v_{n}-
     %\frac{1}{2}(\delta_0 + \phi_0 + 1) v_0 + b_0\\
  &= \frac{1}{2}(\delta_0 + \phi_0 + 1) (v_{n} - v_0) + b_0.
\end{align}
Finally, \eqref{eq:7} is written as follows in terms of $\phi_n$ and $\delta_n$ using
\eqref{eq:solphi} and \eqref{eq:soldelta}, as
\begin{align}
 %\D &= - \delta_{n-1} u_n w_{n-1} + \delta_n u_{n+1} w_n - v_n^2 + \phi_n u_n w_{n-1} -
 %u_n w_{n-1} - \phi_{n+1} u_{n+1} w_n - u_{n+1} w_n,\\
 %&= (\delta_n - \phi_{n+1} - 1)\kappa_{n+1} -(\delta_{n-1} - \phi_n + 1)\kappa_n - v_n^2
 %, \qquad \kappa_n \equiv u_n w_{n-1},\\
 \D + v_n^2 &= (\delta_0 - \phi_0 - 2(n+1))\kappa_{n+1}
                -(\delta_0 - \phi_0 - 2(n-1))\kappa_n,
\end{align}
with
\begin{align}
 \kappa_n \equiv u_n w_{n-1}.\label{kndef}
\end{align}
Multipliying both sides by $(\delta_0 - \phi_0 -2n)$ as an integrating factor, one can
reduce the above to
\begin{align}\label{kappa-notsummed}
 (\delta_0 - \phi_0 -2n)(\delta_0 - \phi_0 - 2n + 2)\kappa_n
 = (\delta_0 - \phi_0)(\delta_0 - \phi_0 +2)\kappa_0
 + \sum\limits_{k=0}^{n-1}(\D + v_k^2)(\delta_0 - \phi_0 -2k).
\end{align}
The sum over $k$ in \eqref{kappa-notsummed} can be reexpressed\footnote{This is
done by noticing the sum to be telescopic or via the polygamma function of the first
order.} as
\begin{align}\label{sumonkend}
 \sum\limits_{k=0}^{n-1}(\D + v_k^2)(\delta_0 - \phi_0 -2k)
 &= \frac{n(\delta_0 - \phi_0 - n+1)
  (\D(\delta_0 - \phi_0 - 2n + 1)^2 + v_0^2(\delta_0 - \phi_0 - 1)^2)}
         {(\delta_0 - \phi_0 - 2n + 1)^2}.
\end{align}
From \eqref{kappa-notsummed} and \eqref{sumonkend}, recalling that $u_0$ was required to
vanish so that $\kappa_0 = u_0 w_{-1} = 0$, one has
\begin{align}\label{kn}
 \kappa_n = \frac{n(\delta_0 - \phi_0 - n+1)(\D(\delta_0 - \phi_0 - 2n + 1)^2
  + v_0^2(\delta_0 - \phi_0 - 1)^2)}
  {(\delta_0 - \phi_0 - 2n + 1)^2(\delta_0 - \phi_0 -2n)(\delta_0 - \phi_0 - 2n + 2)}.
\end{align}

\subsection{The linear pencil $\X+\mu\Z$}
The algebra $\A$ is invariant under the affine transformation
\begin{align*}
 \X \longmapsto \X + \mu\Z, \qquad \mu\in\R.
\end{align*}
As a result, one expects the transformed solutions for the coefficients in
\eqref{eq:actions} to be given by \eqref{eq:solphi}, \eqref{eq:soldelta} and
\eqref{eq:solb} with modified parameters. Indeed one finds the parameters to be replaced
by
\begin{align*}
 \phi_0 \longmapsto \phi_0 - \mu, \qquad
 \delta_0 \longmapsto \delta_0 - \mu, \qquad
 b_0 \longmapsto b_0 - \mu v_0.
\end{align*}
Thus, the diagonalization of the linear pencil $\X+\mu\Z$ amounts to the diagonalization
of $\X$ up to a shift in the parameters.

\subsection{Representations on polynomials}
Denoting by $\bra{x}$ the dual eigenvectors:
\begin{align*}
\bra{x}X = x\bra{x},
\end{align*}
one can look for the polynomials $q_n(x)\equiv\langle x \ket{n}$ that diagonalize $X$
\begin{align}
 X q_n(x) \equiv x\, q_n(x) = c_n q_{n-1}(x)+b_n q_n(x)+a_n q_{n+1}(x).
\end{align}
% Multiplying the above by $\prod_{i=0}^{n-1}a_i$ and thus adopting the
By appropriate renormalization, one obtains a monic recurrence relation
\begin{align}\label{eq:eigen}
 X p_n(x) \equiv x p_n(x) = a_{n-1} c_n p_{n-1}(x)+b_n p_n(x)+p_{n+1}(x), \qquad
 p_n(x)=\left(\prod\limits_{i=0}^{n-1}a_i\right) q_n(x).
\end{align}
% \begin{remark}
The families of polynomials $p_n(x)$ that diagonalize $X$ can be determined by
identifying the coefficients $a_{n-1}c_n$ and $b_n$.
% \end{remark}

From \eqref{eq:solphi}, \eqref{eq:soldelta}, \eqref{kndef} and \eqref{kn}, one
has that
\begin{align}\label{eq:ancn}
 %a_n &= u_n \phi_n, \qquad c_n = w_n \delta_n \quad \implies \quad a_n c_{n-1} = \kappa_n
 %(\phi_0 + n) (\delta_0 - n+1),\\ a_n c_{n-1} &= (\phi_0 + n) (\delta_0 - n+1)
 %\frac{n(\delta_0 - \phi_0 - n+1)(\D(\delta_0 - \phi_0 - 2n + 1)^2 - v_0^2(\delta_0 -
 %\phi_0 - 1)^2)}{(\delta_0 - \phi_0 - 2n + 1)^2(\delta_0 - \phi_0 -2n)(\delta_0 - \phi_0
 %- 2n + 2)}\\
 a_{n-1} c_n &= (n + \phi_0) (n - \delta_0 - 1)
 \frac{n(n + \phi_0 - \delta_0 - 1)
  (\D(2n + \phi_0 - \delta_0 - 1)^2 + v_0^2(\phi_0 - \delta_0 + 1)^2)}
 {(2n + \phi_0 - \delta_0 - 1)^2(2n + \phi_0 - \delta_0)(2n + \phi_0 -\delta_0 - 2)}
\end{align}
and from \eqref{eq:solb} and \eqref{eq:solv}, that
\begin{align}\label{eq:bn}
 b_n = \frac{1}{2} \frac{(\delta_0+\phi_0+1)(\phi_0-\delta_0+1)(\phi_0-\delta_0 - 1) v_0}
 {(2n + \phi_0 - \delta_0 - 1)(2n + \phi_0 -\delta_0 + 1)} + \tilde{b}_0, \qquad
 \tilde{b}_0 \equiv b_0 - \frac{1}{2}(\delta_0 + \phi_0 + 1) v_0.
\end{align}
Finite-dimensional representations of dimension $N+1$ are obtained
% \RED{
if $w_N=0$ since it
follows that $a_N=0$ from \eqref{eq:9}. This implies that $\kappa_{N+1}=0$.
% }
From \eqref{kn}, we see that this is achieved for any value of $\D$ by
\begin{align}\label{eq:trunc1}
 N=(\delta_0 - \phi_0).
\end{align}
If $\D\neq0$, one finds an additionnal pair of solutions given by
\begin{align}\label{eq:trunc2}
 N+1 &=-\frac{1}{2}\left[\phi_0-\delta_0-1\pm(\phi_0-\delta_0+1)v_0\sqrt{-\D^{-1}}\right].
\end{align}

\section{The case $\D=0$: Jacobi polynomials}\label{sec:c0}
With $\D$ vanishing, the coefficient $a_{n-1}c_{n}$ \eqref{eq:ancn} simplifies to
\begin{align}\label{ancn:jacobi}
 a_{n-1} c_n =
 \frac{n(n+\phi_0)(n-\delta_0-1)(n+\phi_0 - \delta_0 -1)(\phi_0 - \delta_0 + 1)^2 v_0^2}
 {(2n +\phi_0 -\delta_0 - 1)^2(2n +\phi_0 - \delta_0)(2n +\phi_0 - \delta_0 - 2)}.
\end{align}
Setting $v_0=2(\phi_0 - \delta_0 + 1)^{-1}$, one identifies the basis vector to be
proportionnal to the Jacobi polynomials $P_n^{(\alpha,\beta)}(x)$ with parameters
\begin{align}\label{eq:jacobiab}
 \alpha = -\delta_0-1, \qquad \beta = \phi_0.
\end{align}
With $\tilde{b}_0=0$, the coefficient $b_n$ of \eqref{eq:bn} is given by
\begin{align}\label{bn:jacobi}
 b_n %&= \frac{1}{2} \frac{(\delta_0 + \phi_0 + 1)(\phi_0 - \delta_0 + 1)(\phi_0 -
     %\delta_0 - 1) v_0}{(2n + \phi_0 - \delta_0 - 1)(2n + \phi_0 -\delta_0 + 1)}\\
  &= \frac{(\beta^2+\alpha^2)}{(2n + \beta + \alpha)(2n + \beta + \alpha +2)}.
\end{align}
Comparing the expressions \eqref{ancn:jacobi} and \eqref{bn:jacobi} for the coefficients
using for instance \cite{KoekoekLeskyetal2010}, we conclude:
\begin{proposition}
In the case $\D=0$, the eigenfunctions $p_{n}(x)$ of $X$ \eqref{eq:eigen} are the
monic Jacobi polynomials
\begin{align*}
 p_n^{(\alpha,\beta)}(x)=\frac{2^n n!}{(n+\alpha+\beta+1)_n}P_n^{(\alpha,\beta)}(x).
\end{align*}
with parameters $\alpha$, $\beta$ given in \eqref{eq:jacobiab}.
\end{proposition}
The only truncation condition possible is \eqref{eq:trunc1}. However, it yields
singular expressions in \eqref{ancn:jacobi} and \eqref{bn:jacobi} for $n\leq N$.

\section{The case $\D>0$: Continuous Hahn polynomials}\label{sec:cpu}
If $\D\neq0$, upon scaling the generators of the algebra according to
\begin{align*}
 \tilde{\X} = \Omega\X, \qquad \tilde{\Z} = \Omega\Z,
\end{align*}
we obtain
\begin{equation}\label{eq:scaling}
 [\tilde{\Z},\tilde{\X}] = \tilde{\Z}^2 + \Omega^{2}\D.
\end{equation}
In view of \eqref{eq:scaling}, one can choose $\Omega$ so that $\D=\pm\tfrac14$.  In this
section, we shall consider the case $\D=+\tfrac14$. The coefficient $a_{n-1}c_n$
\eqref{eq:ancn} is then given by
\begin{align}\label{ancn:conthahn}
 a_{n-1} c_n
 &= (n + \phi_0) (n - \delta_0 - 1)
 \frac{n(n+\phi_0-\delta_0-1)((2n+\phi_0-\delta_0-1)^2/4+v_0^2(\phi_0-\delta_0+1)^2)}
 {(2n + \phi_0 - \delta_0 - 1)^2(2n + \phi_0 - \delta_0)(2n + \phi_0 -\delta_0 - 2)}.
\end{align}
Writing
\begin{align}\label{param-conthahn}
 \phi_0+1 = a+c, \qquad -\delta_0 = b+d, \qquad v_0 = - i\frac{(a-b-c+d)}{2(a+b+c+d)},
\end{align}
one can factorize the term with $v_0$:
\begin{align*}\label{}
 \tfrac14(2n+\phi_0-\delta_0-1)^2+v_0^2(\phi_0-\delta_0+1)^2=(n+a+d-1)(n+b+c-1).
\end{align*}
With \eqref{param-conthahn} and the above, \eqref{ancn:conthahn} becomes
\begin{multline}\label{ancn:conthahngood}
 a_{n-1} c_n = (n + a+c-1) (n +b+d - 1)\\
 \times\frac{n(n + a + b +c +d -2)(n+a+d-1)(n+b+c-1)}
            {(2n + a+b+c+d-1)(2n + a+b+c+d -2)^2(2n + a+b+c+d-3)}.
\end{multline}
Using \eqref{param-conthahn} and taking $\tilde{b}_0=\frac{i}{4}(a+b-c-d)$, the
coefficient $b_n$ \eqref{eq:bn} is found to be
\begin{align}\label{bn:conthahn}
\begin{aligned}
 b_n %&= \frac{1}{2} \frac{(\delta_0 + \phi_0 + 1)(\phi_0 - \delta_0 + 1)(\phi_0 -
% \delta_0 - 1) v_0}{(2n + \phi_0 - \delta_0 - 1)(2n + \phi_0 -\delta_0 + 1)} +
% \tilde{b}_0\\
  %&= -\frac{i}{4} \frac{(a-b+c-d)(a-b-c+d)(a+b+c+d-2)}{(2n + a+b+c+d -2)(2n + a+b+c+d)} +
  %\frac{i}{4}(a+b-c-d)\\
 &= i\Big[-\frac{(n+a+b+c+d-1)(n+a+c)(n+a+d)}{(2n + a+b+c+d -1)(2n + a+b+c+d)}\\
 &\hspace{2em}+\frac{n(n+b+c-1)(n+b+d-1)}{(2n + a+b+c+d -2)(2n + a+b+c+d-1)} + a \Big].
\end{aligned}
\end{align}
The coefficients \eqref{ancn:conthahngood} and \eqref{bn:conthahn} can be identified in
\cite{KoekoekLeskyetal2010} and one arrives at:
\begin{proposition}
In the case $\D>0$, the eigenfunctions $p_{n}(x)$ of $X$ \eqref{eq:eigen} are the
monic continuous Hahn polynomials $P_n^{(a,b,c,d)}(x)$  with parameters given in
\eqref{param-conthahn}.
\end{proposition}

\subsection{Finite-dimensional representations and orthogonal polynomials}
Using \eqref{param-conthahn}, condition \eqref{eq:trunc1} becomes
\begin{align}\label{trunc-parakraw}
 N-1 = -a-c-b-d,
\end{align}
which leads to expressions for \eqref{ancn:conthahn} and \eqref{bn:conthahn} that are
ill-defined for $n<N$. However, this can be resolved using limits and one thus obtains the
para-Krawtchouk polynomials \cite{VinetZhedanov2012}.

Condition \eqref{eq:trunc2} reads
\begin{align}\label{trunc2-chahn}
 N+1= -\frac{1}{2}\left[(a+b+c+d-2)\pm(a-b-c+d)\right] =
 \begin{cases}
  -a-d+1\\
  -b-c+1
 \end{cases}
\end{align}
and corresponds to the truncation of the continuous Hahn polynomials to Hahn
polynomials.

However, for each of these truncations \eqref{trunc-parakraw} and \eqref{trunc2-chahn} to
define real polynomials, the operator $X$ has to be scaled by an imaginary number; this is
equivalent to setting $\D\to-\D$ which corresponds to the situation $\D<0$ that is the
object of the next section.

\section{The case $\D<0$: Hahn and para-Krawtchouk polynomials}\label{sec:cmu}
% In this case, one has
% \begin{align}
%  a_{n-1} c_n &= (n + \phi_0) (n - \delta_0 - 1) \frac{n(n + \phi_0 - \delta_0 -
%  1)(-\rvert\D\rvert(2n + \phi_0 - \delta_0 - 1)^2 + v_0^2(\phi_0 - \delta_0 + 1)^2)}{(2n
%  + \phi_0 - \delta_0 - 1)^2(2n + \phi_0 - \delta_0)(2n + \phi_0 -\delta_0 - 2)}\\ %&=
%  -(n + \phi_0) (n - \delta_0 - 1) \frac{n(n + \phi_0 - \delta_0 -
%  1)(2n-2+(\phi_0-\delta_0+1)(1+2v_0))(2n-2+(\phi_0-\delta_0+1)(1-2v_0))}{4(2n + \phi_0 -
%  \delta_0 - 1)^2(2n + \phi_0 - \delta_0)(2n + \phi_0 -\delta_0 - 2)}
% \end{align}
When $\D<0$, polynomials of a real variable are obtained only if \eqref{eq:trunc1} or
\eqref{eq:trunc2} are satisfied. We begin by treating the latter case.

\subsection{Hahn polynomials}
In view of \eqref{eq:scaling}, we may take $\D=-\tfrac14$ without loss of generality.
Expressing the parameters as follows
\begin{align}\label{eq:param_hahn}
 \phi_0=\beta, \qquad -\delta_0 = \alpha + 1, \qquad
 v_0 = - \frac{(\alpha+\beta+2N+2)}{2(\alpha + \beta +2)}, \qquad
 \tilde{b}_0=\frac{1}{4}(2N -\alpha + \beta),
\end{align}
so that \eqref{eq:trunc2} is satisfied, one obtains
\begin{align}\label{ancn:hahn}
 a_{n-1}c_n&=\frac{n(n+\alpha)(n+\beta)(n+\alpha+\beta)(n+\alpha+\beta+N+1)(N-n+1)}
                  {(2n+\alpha+\beta-1)(2n+\alpha+\beta)^2(2n+\alpha+\beta+1)},
\end{align}
as well as
\begin{align}\label{bn:hahn}
 b_n %&= \frac{1}{4} \frac{(\alpha^2 - \beta^2)(\alpha+\beta+2N+2)}{(2n + \alpha +
 % \beta)(2n + \alpha + \beta + 2)} + \frac{1}{4}(2N -\alpha + \beta)\nonumber\\
 &=\frac{(n+\alpha+\beta+1)(n+\alpha+1)(N-n)}
        {(2n+\alpha+\beta+1)(2n+\alpha+\beta+2)}
  +\frac{n(n+\alpha+\beta+N+1)(n+\beta)}
        {(2n+\alpha+\beta)(2n+\alpha+\beta+1)}.
\end{align}
The coefficients given by \eqref{ancn:hahn} and \eqref{bn:hahn} are found in
\cite{KoekoekLeskyetal2010}.
\begin{proposition}
In the case $\D<0$, the eigenfunctions $p_{n}(x)$ of $X$ \eqref{eq:eigen} related to the
finite-dimensional representation condition \eqref{eq:trunc2} are given in terms of the
monic Hahn polynomials $Q_n^{(\alpha,\beta)}(x)$ for the choice of parameters given in
\eqref{eq:param_hahn}.
\end{proposition}
As previously mentionned, these polynomials can also be obtained as a truncation of
the recurrence defined by \eqref{ancn:conthahn} and \eqref{bn:conthahn}. Indeed, setting
\begin{align}
 \alpha = a+c-1, \qquad \beta = b+d-1,
\end{align}
with one of \eqref{trunc2-chahn}, the coefficient \eqref{ancn:conthahn} and
\eqref{bn:conthahn} become proportional to \eqref{ancn:hahn} and \eqref{bn:hahn},
respectively.
% \begin{align}
%  a_n c_{n-1} %&= 4\frac{n(n + \alpha) (n +\beta)(n + \alpha +
%  \beta)(n+a+d-1)(n+b+c-1)}{(2n + \alpha + \beta -1)(2n + \alpha + \beta )^2(2n +
%  \alpha+\beta +1)}\\ &= -4\frac{n(n + \alpha) (n +\beta)(n + \alpha +
%  \beta)(N-b+1)(n+\alpha+\beta+N+1)}{(2n + \alpha + \beta -1)(2n + \alpha + \beta )^2(2n
%  + \alpha+\beta +1)},
% \end{align}
% while \eqref{bn:conthahn} becomes
% \begin{align}
%  b_n -2ia &=2i\Big[\frac{(n+\alpha+\beta+1)(n+\alpha+1)(N-n)}{(2n + \alpha+\beta+1)(2n +
%  \alpha + \beta +2)}\\ &\qquad\qquad\qquad + \frac{n(n+\alpha + \beta +N +
%  1)(n+\beta)}{(2n + \alpha + \beta)(2n + \alpha + \beta + 1)}\Big].
% \end{align}
%Normalizing the basis vectors by $\ket{n} \to (2i)^{n}\ket{n}$, the action of $i\X/2$ is
%identified as the recurrence relation of the monic Hahn polynomials. This
%normalization is equivalent to setting $\D=\frac{-1}{4}$
Hence, the action of $iX$ when $\D=+\tfrac14$ also leads to the recurrence relation of the
monic Hahn polynomials.

\subsection{Para-Krawtchouk polynomials}
We shall finally indicate how a family of finite-dimensional representations of $\A$
relates to para-Krawtchouk polynomials.  Consider the condition \eqref{eq:trunc1}.
Although leading to singular expressions for certain values of $n$, well-defined
polynomials are obtained by carefully taking limits. Mindful of \eqref{eq:scaling}, it is
convenient in this case to take $\D=-1$. Let $N=2j+p$ with $j$ an integer and $p=0,1$
depending on the parity of $N$, and set
% \begin{align}
%  a_n c_{n-1} &= (n + \phi_0) (n - \delta_0 - 1) \frac{n(n + \phi_0 - \delta_0 - 1)(-(2n
%  + \phi_0 - \delta_0 - 1)^2 + v_0^2(\phi_0 - \delta_0 + 1)^2)}{(2n + \phi_0 - \delta_0 -
%  1)^2(2n + \phi_0 - \delta_0)(2n + \phi_0 -\delta_0 - 2)}\\ &= -(n + \phi_0) (n -
%  \delta_0 - 1) \frac{n(n + \phi_0 - \delta_0 -
%  1)(2n-2+(\phi_0-\delta_0+1)(1+v_0))(2n-2+(\phi_0-\delta_0+1)(1-v_0))}{(2n + \phi_0 -
%  \delta_0 - 1)^2(2n + \phi_0 - \delta_0)(2n + \phi_0 -\delta_0 - 2)}
% \end{align}
% \begin{align}\label{param-parakraw}
%  \phi_0 + 1 = -j+e_1 t, \quad -\delta_0 = -j+e_2 t+1-p, \quad v_0 =
%  -\frac{(2a+2j-2t+2d-1+p)}{2(-2j+e_1t+e_2t-p+1)}, \quad e_1=e_2=1.
% \end{align}
\begin{align}\label{param-parakraw}
 \phi_0 + 1 = -j+e_1 t, \quad
 -\delta_0 = -j+e_2 t+1-p, \quad
 v_0 = \frac{(\gamma+p-1)}{(-2j+e_1t+e_2t-p+1)}, \quad
 e_1=e_2=1.
\end{align}
The parameters $e_1$ and $e_2$ are chosen equal in order to simplify the expressions. The
more general solutions can be recovered using isospectral deformations
\cite{LemayVinetetal2016, GenestTsujimotoetal2017}. With the above parametrization, it can
be seen that \eqref{eq:trunc1} is verified in the limit where $t\to0$. With the parameters
as in \eqref{param-parakraw}, the coefficient $a_{n-1}c_{n}$ \eqref{eq:ancn} becomes
% \begin{multline}
%  a_n c_{n-1} = -(n -j+t-1) (n -j+t - p) \\ \times \frac{n(n -2j+2t -p
%  -1)(2n-2+(-2j+2t-p+1)(1+2v_0))(2n-2+(-2j+2t-p+1)(1-2v_0))}{4(2n -2j+2t -p -1)^2(2n
%  -2j+2t-p)(2n -2j+2t-p -2)}
% \end{multline}
%$$ (2n-2+(-2j+2t-p+1)(1+2v_0))(2n-2+(-2j+2t-p+1)(1-2v_0)) = 4(n-2j+2t-a-d-p)(n-1+a+d)$$
%$$ (2n-2-2j+2t-p+1+2v_0(-2j+2t-p+1))(2n-2-2j+2t-p+1-2v_0(-2j+2t-p+1))$$
%$$= (2n-2-2j+2t-p+1-(2a+2j-2t+2d-1+p))(2n-2-2j+2t-p+1+(2a+2j-2t+2d-1+p)) $$
%$$= (2n-2-4j+4t-2p+2-2a-2d)(2n-2+2a+2d) $$
%$$= 4(n-2j+2t-a-d-p)(n-1+a+d) $$
% \begin{multline}
%  a_n c_{n-1} = -(n -j+t-1) (n -j+t - p) \\
%  \times \frac{n(n -2j+2t -p -1)(n-2j+2t-a-d-p)(n-1+a+d)}{(2n -2j+2t -p -1)^2(2n
%  -2j+2t-p)(2n -2j+2t-p -2)}
% \end{multline}
\begin{multline}
 a_{n-1} c_n = (n -j+t-1) (n -j+t - p) \\
 \times \frac{n(n -2j+2t -p -1)(N-2n+p+\gamma)(N-2n-p+2-\gamma)}
             {(2n -2j+2t -p -1)^2(2n -2j+2t-p)(2n -2j+2t-p -2)}.
\end{multline}
% \begin{multline}
%  \lim_{t\to0} a_n c_{n-1} = -n(n -N -1)(n-N-a-d)(n-1+a+d)\\
%  \times \lim_{t\to0} \frac{(n -j+t-1) (n -j+t - p)}{(2n -2j+2t -p -1)^2(2n -2j+2t-p)(2n
%  -2j+2t-p -2)}
% \end{multline}
% \begin{multline}
%  \lim_{t\to0} \frac{(n -j+t-1) (n -j+t - p)}{(2n -2j+2t -p -1)^2(2n -2j+2t-p)(2n
%  -2j+2t-p -2)}\\ =
%  \begin{cases}
%   \displaystyle\frac{1}{4(2n -2j -1)^2} & p=0\\
%   \displaystyle\frac{1}{4(2n -2j-1)(2n -2j-3)} & p=1
%  \end{cases}\\
%  =\frac{1}{4(2n -2j-p+p-1)(2n -2j-p-p -1)}=\frac{1}{4(2n -N+p-1)(2n -N-p -1)}
% \end{multline}
% \begin{align}
%  \implies \lim_{t\to0} a_n c_{n-1} = \frac{-n(n -N -1)(n-N-a-d)(n-1+a+d)}{4(2n
%  -N+p-1)(2n -N-p -1)}
% \end{align}
Taking the limit $t\to0$ and treating the cases for $p=0,1$ separately, one finds that the
results can be combined as follows
\begin{align}\label{ancn:parakraw}
 \lim_{t\to0} a_{n-1} c_n = \frac{n(N+1-n)(N-2n+p+\gamma)(N-2n-p+2-\gamma)}
                                 {4(2n -N+p-1)(2n -N-p -1)}.
\end{align}
% \begin{align}
%  b_n = \frac{1}{2} \frac{(\delta_0 + \phi_0 + 1)(\phi_0 - \delta_0 + 1)(\phi_0 -
%  \delta_0 - 1) v_0}{(2n + \phi_0 - \delta_0 - 1)(2n + \phi_0 -\delta_0 + 1)} + b
% \end{align}
For the coefficient $b_n$, setting $\tilde{b}_0 = \frac{1}{2}(N+\gamma-1)$ and
inserting \eqref{param-parakraw} in \eqref{eq:bn}, one finds
% \begin{align}
%  b_n = -\frac{1}{4} \frac{(p-1)(-2j-p+2t-1)(2a+2j-2t+2d-1+p)}{(2n -2j-p+2t
%  -1)(2n-2j-p+2t + 1)} + \frac{1}{4}(1-2a-2d-p).
% \end{align}
\begin{align}
 b_n = \frac{1}{2} \frac{(p-1)(-2j-p+2t-1)(\gamma+p-1)}{(2n -2j-p+2t -1)(2n-2j-p+2t + 1)}
     + \frac{1}{2}(N+\gamma-1).
\end{align}
Treating the cases $p=0$ or $p=1$ separately and taking the limit $t\to0$, one sees
that the results can be written jointly as
\begin{align}\label{bn:parakraw}
 \lim_{t\to0} b_n %&=
%  \begin{cases}
% \displaystyle\frac{1}{2} \frac{(N + 1) (\gamma-1)}{(2n -N -1)(2n-N + 1)} +
% \frac{1}{2}(N+\gamma-1), & p=0\\
% \displaystyle\frac{1}{2}(N+\gamma-1) & p=1
%  \end{cases}\nonumber\\
 &= -\frac{(N-n)(N-2n-2+p+\gamma)}{2(2n-N-p+1)}-\frac{n(N-2n+2-p-\gamma)}{2(2n-N+p-1)}.
\end{align}
% \begin{align}
%  \lim_{t\to0} b_n &=
%  \begin{cases}
% \displaystyle\frac{1}{4} \frac{(N - 1) (N+2a+2d-1)}{(2n -N -1)(2n-N + 1)} +
% \frac{1}{4}(1-2a-2d), & p=0\\
% \displaystyle-\frac{1}{2}(a+d) & p=1
%  \end{cases}\\
%  &= -\frac{(n-N)(n+a+d)}{2(2n-N-p+1)}+\frac{n(n-N-a-d)}{2(2n-N+p-1)}
% \end{align}
The coefficients given by \eqref{ancn:parakraw} and \eqref{bn:parakraw} are recognized
in \cite{LemayVinetetal2016} as the coefficients for the recurrence relation of the
monic para-Krawtchouk polynomials.
\begin{proposition}
In the case $\D<0$, the eigenfunctions $p_{n}(x)$ of $X$ \eqref{eq:eigen} in the
finite-dimensional representation \eqref{eq:trunc1} of $\A$ are the monic para-Krawtchouk
polynomials.
\end{proposition}

\section{Conclusion}\label{sec:conclusion}
We have studied tridiagonal representations of the algebra $\A$ with defining relation
$[\Z,\X]=\Z^2+\D$.  Depending on the value of $\D$, in these representations, the linear
pencil $X+\mu Z$ entailed the recurrence relations of the Jacobi ($\D=0$), continuous Hahn
($\D>0$), Hahn and para-Krawtchouk ($\D<0$) polynomials.

In the wake of this work, two research avenues present themselves. One is the exploration
of the tridiagonal representations of the algebra $[\Z,\X]=\Z^{2}+\alpha\X$, another class
of the general quadratic algebra \eqref{eq:genquad}. It is expected that the tridiagonal
representations will lead to the Wilson, Racah and para-Racah polynomials in a similar
fashion.

Another related direction is the study of the so-called meta algebras, poised to describe
both polynomial and rational functions of a given type., as shown in
\cite{VinetZhedanov2020a} for functions of the Hahn type.  The meta-Hahn algebra is in
fact obtained by adjoining to $\A$ an additional generator.  As it turns out, the extra
generator offers a rationale for considering tridiagonal representations.
% This was introduced for the Hahn family in where a single algebra, named the meta-Hahn
% algebra, was seen to encode properties of both the polynomial and rational functions of
% Hahn type.  This meta-Hahn algebra is obtained by These meta algebras are
% three-generated and the addition of a generator
This suggests in particular that the work on the $q$-oscillator algebra
\cite{TsujimotoVinetetal2017} should be revisited in order to bring to the fore the
associated rational functions.

\subsection*{Acknowledgments}
While part of this research was conducted, G. Bergeron and J. Gaboriaud held a scholarship
from the Institut des Sciences Mathématiques of Montreal (ISM) and J. Gaboriaud was partly
funded by an Alexander-Graham-Bell scholarship from the Natural Sciences and Engineering
Research Council of Canada (NSERC).
A. Beaudoin held an Undergraduate Student Research Award (USRA) from the NSERC.
A. Brillant benefitted from a CRM-ISM intern scholarship.
The research of L. Vinet is supported in part by a Discovery Grant from NSERC.
A. Zhedanov who is funded by the National Foundation of China (Grant No.11771015)
gratefully acknowledges the hospitality of the CRM over an extended period and the award
of a Simons CRM professorship.

% \section*{Data availability}
% The data that supports the findings of this study are available within the article.
% \newpage
% \bibliographystyle{abbrv}
% \bibliography{/home/julien/.local/data/all.bib}

% === bibliography
% \newpage
\printbibliography

\end{document}